\documentclass[submission]{pmsart}[2006/01/13]

\usepackage{amsfonts}
\usepackage{amsmath}

\newcommand{\PP}{{\mathbb P}}
\newcommand{\EE}{{\mathbb E}}
\newcommand{\convdistr}{\stackrel{\mathcal{D}}{\rightarrow}}
\newcommand{\convprob}{\stackrel{\mathcal{P}}{\rightarrow}}
\newcommand{\convas}{\stackrel{{\rm a.s.}}{\rightarrow}}

\newcommand{\bearno}{\begin{eqnarray*}}
\newcommand{\enarno}{\end{eqnarray*}}

\volume{0}%
\fasc{0}%
\years{0000}%
\pages{000--000}%
\firstpage{1} %
\lastrevision{xx.xx.xxxx}

\title{A note on first passage probabilities of a L\'{e}vy process reflected at a general barrier}
\titlethanks{}
\dedicated{Dedicated to prof. T. Rolski}
\ShortTitle{L\'{e}vy process reflected at a general barrier} 
\ShortAuthors{Z. Palmowki and P. \'Swi\c{a}tek}
\NumberOfAuthors{2}

\received{31.08.2016}
\FirstAuthor{Zbigniew Palmowski}
\FirstAuthorCity{Wroc\l aw}
\FirstAuthorAffiliation{Faculty of Pure and Applied Mathematics,
Wroc\l aw University of Science and Technology}
\FirstAuthorAddress{Wyb. Wyspia\'nskiego 27, 50-370 Wroc\l aw, Poland}
\FirstAuthorEmail{zbigniew.palmowski@gmail.com}                              
\FirstAuthorThanks{This work is partially supported by the National Science Centre under the grant  DEC-2013/09/B/ST1/01778.
First author kindly acknowledges partial support by the project RARE -318984, a Marie Curie IRSES Fellowship within the 7th European
Community Framework Programme.}                             

\SecondAuthor{Przemys\l aw \'Swi\c{a}tek} 
\SecondAuthorCity{Wroc\l aw} %
\SecondAuthorAffiliation{Credit Suisse}%
\SecondAuthorAddress{pl. Grunwaldzki 25, 50-365 Wroc\l aw, Poland}%
\SecondAuthorEmail{przemek.swiatek@gmail.com}%
\SecondAuthorThanks{}%

\ThirdAuthor{}
\ThirdAuthorCity{}%
\ThirdAuthorAffiliation{}%
\ThirdAuthorAddress{}%
\ThirdAuthorEmail{}
\ThirdAuthorThanks{}%

\FourthAuthor{}
\FourthAuthorCity{}%
\FourthAuthorAffiliation{}%
\FourthAuthorAddress{}%
\FourthAuthorEmail{}%
\FourthAuthorThanks{}%

\FifthAuthor{} 
\FifthAuthorCity{} %
\FifthAuthorAffiliation{}%
\FifthAuthorAddress{} %
\FifthAuthorEmail{}%
\FifthAuthorThanks{}%

\keywords{reflected L\'evy process, Central Limit Theorem, asymptotics, first passage probability}

\MSCcodes{Primary: 60G51, 60G50; Secondary: 60K25.}

\begin{document}
\maketitle
\begin{abstract}
In this paper we analyze a L\'evy process reflected at a general (possibly random) barrier.
For this process we prove Central Limit Theorem for the first passage time. We also
give the finite-time first passage probability asymptotics.

\end{abstract}

\section{Introduction}\label{sec:intro}
In this paper we consider general L\'evy process reflected at a general (possibly random) barrier.
In dimensions larger than one the existence of the reflected process (even for Brownian motion in a time-independent
domain) is a non-trivial matter. Some smoothness assumptions on the domain are required (see \cite{Tanaka}).
However, in one-dimension there are no such problems (see \cite{Kellathoughts}).

The motivation for considering reflected L\'evy process
comes mainly from various storage, queueing and telecommunication models (see \cite{asmbook, abook}).
Usually one deals with reflecting at the constant zero barrier producing the workload of the L\'evy-driven queue (see \cite{krzysiekmandjes} for survey).
There are many applications where more general reflected L\'{e}vy processes appear though.
One can consider for example models
where a constant input is added to the netput mechanism.
This additional a constant input is not available on liquid basis
and can be removed only after some (possibly random) maturity date.
Therefore it can be described in the model by some downward barrier (see e.g. \cite{KOM}).

There has been a great deal of work on overflow probabilities
which correspond to the first passage time probability of the reflected L\'evy process.
A corresponding interpretation for the Cram\'er-Lundberg risk model is also
possible, where the general barrier represents the future threshold for (continuous) dividend pay out.
Another motivation for considering the reflected L\'evy process comes from
the structural biology which is treated thoroughly in \cite{8}.
Finally, considering the reflection of a L\'evy process provides an example of a time
inhomogeneous process where we can derive quite precise asymptotic results of the first passage probabilities.

In this paper we consider general reflected L\'evy process $X_t$ from another, general and independent L\'evy process $Y_t$
and reflected L\'evy process $X_t$ from the deterministic barrier.
We denote this process by $V_t$.

The reflection from the independent L\'evy processes has been already considered in some particular cases.
For example Burdzy and Nualart \cite{5} were interested in the reflection of a one-dimensional
Brownian motion at the path of another Brownian motion (see also Section 3 in \cite{6}
for general issues regarding the reflection of one-dimensional Brownian motion).
The reflected L\'evy process from the deterministic barrier was also already considered.
It was done in the seminal paper of Hansen \cite{Hansen}.
He derived the Cram\'er asymptotics for the first passage probability, where
deterministic, time-dependent barrier tends to $-\infty$ (see also \cite{10} and \cite{7}).

We continue this research deriving (in both cases
- reflecting from the random process and reflecting from the deterministic barrier)
the Central Limit Theorem for the first passage time
$$\tau(u)=\inf\{t\geq 0: V_t>u\}.$$
We also identify the tail-asymptotics of its distribution, where the time horizon goes
to infinity depending linearly on the crossing level.

The most intriguing and the main result of this paper is given in Theorem \ref{LDP}.
In this theorem we consider a deterministic, increasing and convex barrier $f$
and L\'evy process $X$ being compound Poisson process with negative drift:
$$X_t=\sum_{i=1}^{N_t}E_i-t.$$
Above, all jumps $E_i$ have a distribution function $F$
with the following asymptotics
$$\lim_{x\to\infty}\frac{1-F(x)}{l(x)x^{-\beta}e^{-\alpha x}}=1
$$
for some slowly varying function $l$, $\beta\in \mathbb{R}$ and $\alpha >0$.
The process $N_t$ is an independent of $\{E_i\}$ Poisson process with intensity $\lambda$.
In this case
$$V_t=X_t-\min(0,\inf_{0\leq s\leq t}(X_s-f(s))).$$
When $X_t-f(t)\convas -\infty$ as $t\rightarrow\infty$, under some other technical assumptions,
we find the following logarithmic asymptotics:
$$
\lim_{u\to\infty}\frac{1}{u}\log \PP\left(\tau(u)\leq f^{-1}(cu)\right)=-\alpha (1-c),
$$
where $f^{-1}$ is a right-inverse of the lower barrier $f$ and $0<c<1$.
The proof of the above statement is based on identifying precise asymptotic lower and upper bounds
by constructing simpler, 'envelope-like' processes of $V_t$ and analyzing their first passage-times.
We use some geometrical observation and the fluctuation theory of L\'evy processes in the proof.
The proof is surprising more complex and we do not know how to transfer it to more general L\'evy processes.

We organize the paper in the following way. Firstly, in Section \ref{sec:CLT}, we demonstrate the results concerning Central Limit Theorem for both cases
(reflecting from random and deterministic barriers). Then,  in Section \ref{sec:asymptotyka}, we focus on the asymptotics of the first-passage time
probabilities and prove the main result of Theorem \ref{LDP}.

\section{Central Limit Theorem}\label{sec:CLT}
\subsection{Reflection from the independent L\'evy process}
On $(\Omega, \mathcal{F}, \{\mathcal{F}_{t}\}_{\{t\geq 0\}}, \mathbb{P}) $ we define
two independent on each other general L\'{e}vy processes $\{X_t\}_{t\geq 0}$
and $\{Y_t\}_{t\geq 0}$. The filtration $\mathcal{F}_{t}$ is a natural right-continuous
filtration generated by $(X_t, Y_t)$ satisfying usual conditions.
We assume in this paper that $\mathbb{E}X_1=m_X>0$ and that $\mathbb{E}Y_1=m_Y<\infty$.
We denote $Z_t=X_t-Y_t$. Observe that $Z$ is also a L\'{e}vy process.
We define the reflected process $X$ from the barrier $Y$ in the following way:
\begin{equation} \label{odbity2}
V_t=X_t-\inf_{0\leq s\leq t}(X_s-Y_s)=Z_t-\inf_{0\leq s\leq t}(Z_s)+Y_t.
\end{equation}

Our main interest in this section is the Central Limit Theorem for the first passage time:
\begin{equation}\label{taudef}\tau(u)=\inf\{t\geq 0: V_t>u\}.\end{equation}
By $\convdistr$, $\convprob$ and $\convas$ we denote the convergence in distribution, in probability and almost sure, respectively.
\begin{theorem}\label{tw2}
We have:
\begin{enumerate}
\item[i)] If $\mathbb{E}Z_1>0$ then $$\frac{\tau(u)}{u}\convprob \frac{1}{m_X},$$ $$\frac{\tau(u)-\frac{u}{m_X}}{\sqrt{u}}\convdistr N\left(0,\frac{\omega_X^2}{m_X^3}\right)$$
as $u\rightarrow \infty$, where $\omega_Y^2={\rm Var} (Y_1)$
and $N(m\, \sigma^2)$ denotes a Gaussian random variable with mean $m$ and variance $\sigma^2$.
\item[ii)] If $\mathbb{E}Z_1<0$ then $$\frac{\tau(u)}{u}\convprob \frac{1}{m_Y},$$ $$\frac{\tau(u)-\frac{u}{m_Y}}{\sqrt{u}}\convdistr N\left(0,\frac{\omega_Y^2}{m_Y^3}\right)$$
as $u\rightarrow \infty$.
\end{enumerate}
\end{theorem}

In the proof we will need the following lemma.
\begin{lemma}\label{tw1}
\begin{enumerate}
\item[i)] If $\mathbb{E}Z_1>0$ then $$\frac{V_t-m_Xt}{\sqrt{t}}\convdistr N(0,\omega_X^2)$$
as $t\rightarrow \infty$.
\item[ii)] If $\mathbb{E}Z_1<0$ then $$\frac{V_t-m_Yt}{\sqrt{t}}\convdistr N(0,\omega_Y^2)$$
as $t\rightarrow \infty$.
\end{enumerate}
\end{lemma}
\texttt{Proof}. If $\mathbb{E}Z_1>0$ then $\inf_{0\leq s\leq t}(Z_s)\convdistr Z$ as $t\rightarrow \infty$ for some finite random variable $Z$. Hence:
$$\frac{V_t-m_Xt}{\sqrt{t}}=\frac{X_t-\inf_{0\leq s\leq t}Z_s-m_Xt}{\sqrt{t}}=\frac{X_t-m_Xt}{\sqrt{t}}-\frac{\inf_{0\leq s\leq t}Z_s}{\sqrt{t}}\convdistr N(0,\omega_X^2).$$
On the other hand, if $\mathbb{E}Z_1<0$ then $Z_t-\inf_{0\leq s\leq t}Z_s\convdistr Z'$ as $t\rightarrow \infty$ for some finite random variable $Z'$. Thus
\begin{eqnarray*}
\frac{V_t-m_Yt}{\sqrt{t}}&=&\frac{Z_t-\inf_{0\leq s\leq t}Z_s+Y_t-m_Yt}{\sqrt{t}}\\&=&\frac{Z_t-\inf_{0\leq s\leq t}Z_s}{\sqrt{t}}+\frac{Y_t-m_Yt}{\sqrt{t}}\convdistr N(0,\omega_Y^2)\end{eqnarray*}
as $t\rightarrow \infty$.
\hfill $\square$

\texttt{Proof of Thm. \ref{tw2}}.
Let $\mathbb{E}Z_1>0$. Recall that $\inf_{0\leq s\leq t}(Z_s)\convdistr Z$ as $t\rightarrow \infty$ for some finite random variable $Z$.
From the Law of Large Numbers we have:
$$m_X=\lim_{t\to\infty}\frac{X_t}{t}=\lim_{t\to\infty}\frac{V_t+\inf_{0\leq s\leq t}Z_s}{t}=\lim_{u\to\infty}\frac{V_{\tau(u)}}{\tau(u)}=\lim_{u\to\infty}\frac{u}{\tau(u)}.$$
Now from Lemma \ref{tw1} it follows that:
$$\frac{V_t-m_Xt}{\sqrt{t}}\convdistr N(0,\omega_X^2)$$
and from Anscombe's Theorem (see e.g. \cite[Thm. 7.3.2]{246}) we get that
$$\frac{V_{\tau(u)}-m_X\tau(u)}{\sqrt{\tau(u)}}\convdistr N(0,\omega_X^2)\; {\rm as}\; u\to\infty.$$
Thus
$$\frac{\sqrt{u}}{\sqrt{\tau(u)}}\frac{u-m_X\tau(u)}{\sqrt{u}}\convdistr N(0,\omega_X^2)$$
and finally
$$\frac{\frac{u}{m_X}-\tau(u)}{\sqrt{u}}\convdistr N\left(0,\frac{\omega_X^2}{m_X^3}\right).$$
The proof of the case $ii)$ is similar.
\hfill $\square$

\subsection{Reflection from the deterministic barrier}

Our goal is to get similar results for a L\'{e}vy process $X$ reflected at some fixed, deterministic, time-dependent and c\'{a}dl\'{a}g boundary $f(t)$.

By $V_t$ we denote the general L\'evy process $X_t$ reflected at the barrier $f(t)$, that is,
\begin{equation} \label{odbity1}
V_t=X_t-\min(0,\inf_{0\leq s\leq t}(X_s-f(s))).
\end{equation}
Thus for $f(0)\geq 0$ we have:
$$V_t=X_t-\inf_{0\leq s\leq t}(X_s-f(s)).$$

For such process we can prove the Central Limit Theorem for $\tau(u)$ defined in (\ref{taudef}).
\begin{lemma}\label{tw3}
\begin{enumerate}
\item[i)] If $X_t-f(t)\convas +\infty$, then $$\frac{V_t-m_Xt}{\sqrt{t}}\convdistr N(0,\omega_X^2).$$
\item[ii)] If $X_t-f(t)\convas -\infty$, then $$V_t\sim f(t),$$
where we write
\begin{equation}\label{sim}
h\sim g\quad {\rm if}\quad \lim_{t\to\infty}h(t)/g(t)=1.
\end{equation}
\end{enumerate}
\end{lemma}
The proof of above lemma is the same as the proof of Lemma \ref{tw1}.

The main result of this section is the following theorem.
\begin{theorem}\label{tw4}
\begin{enumerate}
\item[i)] If $X_t-f(t)\convas +\infty$ as $t\rightarrow \infty$ then $$\frac{\tau(u)}{u}\convprob \frac{1}{m_X},$$ $$\frac{\tau(u)-\frac{u}{m_X}}{\sqrt{u}}\convdistr N\left(0,\frac{\omega_X^2}{m_X^3}\right)$$
as $u\rightarrow \infty$.
\item[ii)] If $X_t-f(t)\convas -\infty$ as $t\rightarrow \infty$
then $$\frac{\tau(u)}{f^{-1}(u)}\convprob 1$$
as $u\rightarrow \infty$ where $f^{-1}(u)=\inf\{t\geq 0: f(t)\geq u\}$.
\end{enumerate}
\end{theorem}

Let $$\xi(u)=X_{\tau(u)}-u.$$
The proof of above theorem is based on the following fact.
\begin{lemma}\label{tw4b}
If $\EE X_1>0$ and $X_t-f(t)\to+\infty$ a.s. as $t\rightarrow \infty$ then $\tau(u)$ and $\xi(u)$ are asymptotically independent:
for bounded and continuous function $h$ on $[0,\infty)$ and bounded continuous function $g_u$ (possibly dependent on $u$) on real line
for $u>0$, we have that:
\begin{equation}\label{stam}
\EE\left(h(\xi(u))g_u\left(\frac{\tau(u)-\frac{u}{m_X}}{\sqrt{u}}\right)\right)\sim\EE f(\xi(\infty))\EE g_u(N),
\end{equation}
where $N$ is a $N\left(0,\frac{\omega^2_X}{m^3_X}\right)$ random variable.
\end{lemma}

\texttt{Proof}. The proof of this lemma is the same as the proof of \cite[Lemma 5.8, p. 368]{asmbook}.
In particular, note that this proof does not change if one allows the dependence of the function $g$ on $u$.
\hfill $\square$

\texttt{Proof of Thm \ref{tw4}.} The case (i) follows straightforward from Lemma \ref{tw4b}.
The case (ii) is a consequence of the observation that $\mathbb{P}(\tau\leq f^{-1}(u))\rightarrow 1$ as $u\rightarrow\infty$.
\hfill $\square$

\begin{remark}\rm
When $X_t-f(t)\convas -\infty$ then $\xi(\infty)$ may not exists.
Indeed, if we take $f(t)=[t]^2$ (where $[t]$ denotes the integer part of $t$) and take $X_t$ being a Poisson process with rate $\lambda$,
then we have:
$$\lim_{n\to\infty}\xi\left(n^2-\frac{1}{2}\right)\convprob \frac{1}{2},$$
$$\lim_{n\to\infty}\xi\left(n^2+\frac{1}{2}\right)\convprob A,$$
where $\PP(A=\infty)=e^{-\lambda}$ and $\PP(A=\frac{1}{2})=1-e^{-\lambda}$.
\end{remark}

\section{Finite-time first passage probabilities: asymptotics}\label{sec:asymptotyka}

\subsection{Reflection from the independent L\'evy process}
We start our considerations, as before, from considering the first passage probabilities
for the general L\'evy process $X$ reflected from another, independent L\'evy process $Y$; see (\ref{odbity2})
for the formal definition.
In this case we can derive the following result.

\begin{theorem}
\begin{enumerate}
\item[i)]
If $\EE Z_1>0$ then, as $u\rightarrow\infty$, for general function $g$ we have
$$\PP(\tau(u)\leq g(u))=\PP\left(\frac{\tau(u)-\frac{u}{m_X}}{\sqrt{u\frac{\omega_X^2}{m_X^3}}}\leq\frac{g(u)-\frac{u}{m_X}}{\sqrt{u\frac{\omega_X^2}{m_X^3}}}\right)\sim \mathbb{P}\left(N(0,1)\leq\frac{g(u)-\frac{u}{m_X}}{\sqrt{u\frac{\omega_X^2}{m_X^3}}}\right),$$
where symbol $\sim$ is defined in (\ref{sim}) and $N(0,1)$ is a standard Gaussian random variable.
\item[ii)]
Similarly, if $\EE Z_1<0$ and $u\rightarrow\infty$ then
$$\PP(\tau(u)\leq g(u))=\PP\left(\frac{\tau(u)-\frac{u}{m_Y}}{\sqrt{u\frac{\omega_Y^2}{m_Y^3}}}\leq\frac{g(u)-\frac{u}{m_Y}}{\sqrt{u\frac{\omega_Y^2}{m_Y^3}}}\right)\sim \mathbb{P}\left(N(0,1)\leq\frac{g(u)-\frac{u}{m_Y}}{\sqrt{u\frac{\omega_Y^2}{m_Y^3}}}\right).$$
\end{enumerate}
\end{theorem}
\texttt{Proof}.
The assertion of the theorem is a consequence of Theorem \ref{tw2}.
\hfill $\square$
\\[1cm]

\subsection{Reflection from the deterministic barrier $f$}
Using Theorem \ref{tw4} similar results could be derived for the deterministic boundary.

\begin{theorem}
Assume that
\begin{equation}\label{as1}
X_t-f(t)\convas \infty
\end{equation}
and that $g$ is a function such that $g(u)\leq f^{-1}(u)$.
Then we have:
$$\PP(\tau(u)\leq g(u))=\PP\left(\frac{\tau(u)-\frac{u}{m_X}}{\sqrt{u\frac{\omega_X^2}{m_X^3}}}\leq\frac{g(u)-\frac{u}{m_X}}{\sqrt{u\frac{\omega_X^2}{m_X^3}}}\right)\sim N\left(\frac{g(u)-\frac{u}{m_X}}{\sqrt{u\frac{\omega_X^2}{m_X^3}}}\right).$$
\end{theorem}

From now on we assume complementary condition that
\begin{equation}\label{as2}
X_t-f(t)\convas -\infty.
\end{equation}

Moreover, we assume that
\begin{equation}\label{convex}
f(t)\quad \text{ is an increasing and convex function;}
\end{equation}
and that $X$ is a compound Poisson process given by:
\begin{equation}\label{comp}
X_t=\sum_{i=1}^{N_t}E_i-t,\end{equation}
where generic $E$ has a distribution function $F$ and $N_t$ is independent of $\{E_i,\ i \in \mathbf{N}\}$
Poisson process with intensity $\lambda$.
Additionally, we will assume that
\begin{equation}\label{F}
\overline{F}(x)=1-F(x)\sim l(x)x^{-\beta}e^{-\alpha x}
\end{equation}
for some slowly varying function $l$ and $\beta\in \mathbb{R}$, $\alpha >0$.

From \cite{BGT} we know that
the slowly varying function $l$ has the following representation:
\begin{equation}\label{slowlyrepr}l(x)=\gamma(x)e^{\int_a^x\frac{h(u)}{u}du},\qquad x\geq a
\end{equation}
for some $a>0$, where $\gamma(\cdot)$ is measurable and $\gamma(x)\to \gamma\in (0,\infty)$, $h(x)\to 0$ as $x\to 0$.

We also assume that:
\begin{equation}\label{xgammaprime}
x\gamma '(x)\to 0\quad \text{when}\quad x\to\infty.
\end{equation}

\begin{remark}\rm
Note that for example when we take $\gamma(x)=1$, $a=e$ and $h(x)=\frac{1}{\log(x)}$ then we will get $l(x)=\log (x)$.

This case is the most difficult one since the shape of the barrier has
influence on the behaviour of the reflected process
(please compare e.g. with \cite{Hansen} where (\ref{as1}) is assumed).
Therefore we manage to derive the asymptotics only under above
restricted assumptions. We believe though that at least logarithmic asymptotic holds true for more general setup.
\end{remark}

The next theorem is the main result of this paper.
\begin{theorem}\label{LDP}
Assume that (\ref{as2})-
(\ref{xgammaprime}) hold and $\beta > 3$. Then for any $c \in (0,1)$:
$$
\lim_{u\to\infty}\frac{1}{u}\log \PP\left(\tau(u)\leq f^{-1}(cu)\right)=-\alpha (1-c).
$$
\end{theorem}
\begin{remark}
\rm In fact in the proof we derive more exact lower and upper uniform bounds which seems to be of own interest.
\end{remark}
We divide the proof of the theorem into two parts: proving lower and upper bounds for
$\PP\left(\tau(u)\leq f^{-1}(cu)\right)$, both of which hold without the assumption of $\beta >3$.

{\bf Lower bound.}
\begin{proposition}\label{lowerbound}
$$\lim_{u\to\infty} \frac{\PP(\tau(u)\leq f^{-1}(cu))}{\lambda e^{-\alpha(1-c) u}\frac{l(u-cu)}{\alpha f'(f^{-1}(cu))}(u-cu)^{-\beta}}\geq 1.$$
\end{proposition}
\texttt{Proof}.
To get the lower bound of $\PP(\tau(u)\leq f^{-1}(cu))$ we consider the following process:
$$\underline{V}_t=f(t)+\Delta V_t,$$
where $\Delta V_t$ is the size of the jump of $V_t$ that occurs exactly at time $t>0$.
Observe that $V_t\geq \underline{V}_t$.
Indeed, note that $V_t-V_{t-}+\Delta V_t$ and $V_{t-}\geq f(t)$ by the definition of the reflected proces $V$.
Denote now $\underline{\tau}(u)=\inf\{t\geq 0: \underline{V}_t>u\}.$
It is clear that $$\PP(\tau(u)\leq f^{-1}(cu))\geq \PP(\underline{\tau}(u)\leq f^{-1}(cu)).$$
We also introduce events $A_k$ that there are exactly $k$ jumps of $\underline{V}_t$ (and also $V_t$) before time $f^{-1}(cu)$.
Then:
$$\PP(\underline{\tau}(u)\leq f^{-1}(cu))=\sum_{k=1}^\infty\PP(\underline{\tau}(u)\leq f^{-1}(cu)|A_k)\PP(A_k),$$
where $\PP(A_k)=\frac{(\lambda f^{-1}(cu))^k}{k!}e^{-\lambda f^{-1}(cu)}$. Besides that conditioned on $A_k$ the times of jumps
have the same law as the order statistics of $k$ independent uniform distributed random variables on the interval $[0,f^{-1}(cu)]$.
Thus:
\begin{eqnarray*}
\PP(\underline{\tau}(u)\leq f^{-1}(cu)|A_k)
=1-\left[1-\frac{1}{f^{-1}(cu)}\int_0^{f^{-1}(cu)}\overline{F}(u-f(x))dx\right]^k
\end{eqnarray*}
and
\begin{eqnarray*}
\lefteqn{\PP(\underline{\tau}(u)\leq f^{-1}(cu))}\\ &&= \sum_{k=1}^\infty \frac{(\lambda f^{-1}(cu))^k}{k!}e^{-\lambda f^{-1}(cu)}\left[1-\left[1-\frac{1}{f^{-1}(cu)}\int_0^{f^{-1}(cu)}\overline{F}(u-f(x))dx\right]^k\right]\\
&&=1-e^{-\lambda\int_0^{f^{-1}(cu)}\overline{F}(u-f(x))dx}.
\end{eqnarray*}
Hence:
\begin{equation}\label{asymF}
\PP(\underline{\tau}(u)\leq f^{-1}(cu))\sim \lambda\int_0^{f^{-1}(cu)}\overline{F}(u-f(x))dx. 
\end{equation}
To identify more precisely above asymptotics we need to analyze asymptotic behaviour of
\begin{equation}\label{asymF2}
\lambda \int_0^{f^{-1}(cu)}l(u-f(x))(u-f(x))^{-\beta}e^{-\alpha(u-f(x))}dx.
\end{equation} Using change of variables $f(x)=z$ we get that:
\begin{eqnarray*}
\lefteqn{\int_0^{f^{-1}(cu)}l(u-f(x))(u-f(x))^{-\beta}e^{-\alpha(u-f(x))}dx}\\&& =e^{-\alpha u}\int_{f(0)}^{cu}l(u-z)(u-z)^{-\beta}e^{\alpha z}\frac{1}{f'(f^{-1}(z))}dz.
\end{eqnarray*}
We denote $g(u,z):=l(u-z)(u-z)^{-\beta}e^{\alpha z}\frac{1}{f'(f^{-1}(z))}$ and
\begin{eqnarray*}
I(u)=\int_{f(0)}^{cu}g(u,z)dz,
\end{eqnarray*}

\begin{eqnarray*}
J(u)=\frac{g(u,cu)}{\alpha}.
\end{eqnarray*}
Now using de l'Hospital rule we will show that:
\begin{eqnarray*}
\frac{I(u)}{J(u)}\to 1\qquad \textrm{as}\ u\to\infty.
\end{eqnarray*}
Indeed, note that
\begin{eqnarray*}
\frac{dI(u)}{du}&=&l(u-cu)(u-cu)^{-\beta}e^{\alpha cu}\frac{c}{f'(f^{-1}(cu))} \\ &&-\beta\int_{f(0)}^{cu}l(u-z)(u-z)^{-\beta-1}e^{\alpha z}\frac{1}{f'(f^{-1}(z))}dz \\ &&+\int_{f(0)}^{cu}\frac{\partial l(u-z)}{\partial u}(u-z)^{-\beta}e^{\alpha z}\frac{1}{f'(f^{-1}(z))}dz.
\end{eqnarray*}
Moreover, we have
$$\frac{\partial l(u-z)}{\partial u}=l(u-z)\frac{\partial \gamma(u-z)}{\partial u}\frac{1}{\gamma(u-z)}-l(u-z)\frac{h(u-z)}{u-z}.$$
Hence we can write:
\begin{eqnarray*}
\frac{dI(u)}{du}&=&cg(u,cu)-\int_{f(0)}^{cu}g(u,z)\frac{\beta}{u-z}dz \\ &&+\int_{f(0)}^{cu}g(u,z)\frac{1}{\gamma(u-z)}\frac{\partial \gamma(u-z)}{\partial u}dz-\int_{f(0)}^{cu}g(u,z)\frac{h(u-z)}{u-z}dz.
\end{eqnarray*}
Last two terms are asymptotically negligible because we assumed in (\ref{slowlyrepr}) and (\ref{xgammaprime}) that both $x\gamma '(x)\to 0$ and $h(x)\to 0$ when $x\to\infty$.
Note also that the first increment dominates the second one.
Furthermore, we can observe that
\begin{eqnarray*}
\frac{dJ(u)}{du}&=& -\frac{\beta(1-c)}{\alpha(u-cu)}g(u,cu)+cg(u,cu)-g(u,cu)\frac{cf''(f^{-1}(cu))}{\alpha^2 f'(f^{-1}(cu))^2}\\ &+&\frac{\gamma '(u-cu)(1-c)}{\alpha\gamma(u-cu)}g(u,cu)-(1-c)\frac{h(u-cu)}{\alpha(u-cu)}g(u,cu)
\end{eqnarray*}
is asymptotic equivalent to $cg(u,cu)$ when $\frac{f''(f^{-1}(cu))}{f'(f^{-1}(cu))^2}\to 0$ as $u\to \infty$. We can rewrite this condition as follows:
\begin{eqnarray*}
\frac{f''(z)}{f'(z)^2}=\left[-\frac{1}{f'(z)}\right]'\to 0 \qquad \textrm{as} \qquad z\to \infty.
\end{eqnarray*}
 It is satisfied for any increasing convex function $f(z)$ since $-\frac{1}{f'(z)}$ is negative increasing function bounded by $0$. Thus it must have a limit as $z\to\infty$
 This implies that derivative of this function vanishes at $\infty$.

We proved thus that both $\frac{dI(u)}{du}$ and $\frac{dJ(u)}{du}$ are asymptotic equivalent to $cg(u,cu)$ as $u\to\infty$. We can now conclude that
\begin{eqnarray*}
\PP(\underline{\tau}(u)\leq f^{-1}(cu)) &\sim& \lambda e^{-\alpha u} I(u) \sim \lambda e^{-\alpha u} J(u) =
\\&=& \lambda l(u-cu)(u-cu)^{-\beta}e^{-\alpha((1-c)u)}\frac{1}{\alpha f'(f^{-1}(cu))}.
\end{eqnarray*}
which completes the proof.
\hfill $\square$

{\bf Upper bound.}
\begin{proposition}\label{upperbound}
For any $0<\epsilon <\alpha$ there exists $\theta^\epsilon$ and constant $C_0$
such that
\begin{eqnarray*}
\PP(\tau(u)\leq f^{-1}(cu))\leq \frac{C_0}{(1+u-cu)^2} e^{-(\alpha-\epsilon)(u-cu-f^{-1}(cu))-\theta^\epsilon f^{-1}(cu)}.
\end{eqnarray*}
\end{proposition}
\texttt{Proof}.
To get the upper bound we consider the following process:
$$\overline{V}_t=cu+f^{-1}(cu)+X_t.$$
Note that $\overline{V}$ will never reach barrier $f$ before time $f^{-1}(cu)$ because it starts at the level $cu+f^{-1}(cu)$. This observation leads to the conclusion that:
\begin{eqnarray}
\PP(\tau(u)\leq f^{-1}(cu))&\leq&\PP(\max_{0\leq s\leq f^{-1}(cu)}\overline{V}_s\geq u)\nonumber\\ &=&\PP(\max_{0\leq s\leq f^{-1}(cu)}X_s\geq u-cu-f^{-1}(cu))\nonumber\\
&=&\PP(\overline{\tau}(u-cu-f^{-1}(cu))\leq f^{-1}(cu)),\label{numer2}
\end{eqnarray}
where
$$\overline{\tau}(x)=\inf\{t\geq 0: X_t>x\}.$$
To estimate above probability we will use asymptotics of the renewal function given in \cite{Hoglund}.
We denote by $L$ local time related to running supremum of $X$ and by $L^{-1}$ its right-continuous inverse. Finally, let $H_t=X_{L^{-1}_t}$.
Process $\{(L_t^{-1}, H_t), t\geq 0\}$ is an ascending ladder height process. We introduce also the sequence $\{e_i\}_{i\in \mathbf{N}}$ of independent exponentially distributed random variables
with intensity $q$. Let $\sigma_n=\sum_{i=1}^ne_i$. Now we consider a two-dimensional random walk
$$\{(S_i,R_i)=(L_{\sigma_i}^{-1}, H_{\sigma_i}),i=1,2,...\}$$ starting from $(0,0)$ and having
the step size distribution:
\begin{eqnarray*}
\mu(dx,dt)=\PP(H_{\sigma_1}\in dx,L_{\sigma_1}^{-1}\in dt).
\end{eqnarray*}
Denote:
\begin{eqnarray*}
\Upsilon(x)=\min\{n: S_n>x\},\qquad G(x,t)=\PP(\Upsilon(x)<\infty, R_{\Upsilon(x)}\leq t)
\end{eqnarray*}
Then from \cite{jaMartijn} we have:
\begin{eqnarray}\label{lematZbyszkowy}
\PP(\overline{\tau}(x)\leq t) \leq \frac{G(x,t+M)}{h(0-,M)}
\end{eqnarray}
for any $M,x,t>0$, where
\begin{eqnarray*}
G(x,t)=\sum_{n=0}^\infty \mu^{\star n}\star h(x,t)
\end{eqnarray*}
is a renewal function and $h(x,t):=\PP(H_{\sigma_1}>x,L_{\sigma_1}\leq t)$ with $h(0-,M)=\lim_{x\uparrow 0}h(x,M)$.

We denote by $\Pi_H$ the L\'evy measure associated with the ladder height process $H$.
Then by Vigon's formula (see e.g. \cite[Cor. 7.9]{abook})
for the spectrally positive L\'evy process $X$ we have:
\begin{eqnarray*}
\overline{\Pi}_H(y)=\int_0^\infty e^{-\Phi(0)z}\overline{\Pi}_X(z+y)dz,
\end{eqnarray*}
where
$$\Phi(\theta)=\inf\{\kappa\geq 0: \log \EE e^{-\kappa X_1}\geq \theta\}$$
is the right inverse of the Laplace exponent of the dual process $-X$.
Hence:
\begin{eqnarray}
\nonumber\frac{\overline{\Pi}_H(y)}{\overline{\Pi}_X(y)}=\int_0^\infty e^{-\Phi(0)z}\frac{\overline{\Pi}_X(z+y)dz}{\overline{\Pi}_X(y)}\sim\int_0^\infty e^{-\Phi(0)z}e^{-\alpha z}dz=\frac{1}{\alpha+\Phi(0)}.\\
\label{numer1}
\end{eqnarray}
Moreover, from \cite[Lem. 3.5]{kyprklupmaller}
\begin{eqnarray*}
\PP(H_{\sigma_1}>x)=q\int_0^{\infty}e^{-qt}\PP(H_t>x)dt\sim q\int_0^{\infty}e^{-qt}t(\EE e^{\alpha H_1})^{t-1}\overline{\Pi}_H(x)dt,
\end{eqnarray*}
From (\ref{numer1}) we can conclude now that:
\begin{eqnarray}\label{prawieVigon}
\PP(H_{\sigma_1}>x)=C_1\overline{\Pi}_H(x)\sim C_2\overline{\Pi}_X(x)\sim C_2\lambda l(x)x^{-\beta}e^{-\alpha x},
\end{eqnarray}
where $C_1=q\int_0^{\infty}e^{-qt}t(Ee^{\alpha H_1})^{t-1}dt$ and $C_2=\frac{C_1}{\alpha+\Phi(0)}$.
Denote
\begin{eqnarray*}
\Theta=\{(\theta_1,\theta_2)\in\mathbb{R}^2 : \int|(x,t)|e^{\theta_1x+\theta_2t}\mu(dx,dt)<\infty\}.
\end{eqnarray*}
From the equation (\ref{prawieVigon}) it follows that for any $0<\epsilon<\alpha$ there exists $\theta^\epsilon$ such that $(\alpha-\epsilon,\theta^\epsilon)\in\Theta$
and
\begin{equation*}
\EE e^{(\alpha-\epsilon)H_{\sigma_1}+\theta^\epsilon L^{-1}_{\sigma_1}} =1.\end{equation*}
Moreover, then we have also
$$\EE\left[(H_{\sigma_1}+L^{-1}_{\sigma_1})^2 e^{(\alpha-\epsilon)H_{\sigma_1}+\theta^\epsilon L^{-1}_{\sigma_1}}\right]<\infty.$$
Hence from \cite[Thm. 1.2]{Hoglund} applied with $(c_1,c_2)=(1,1)$
we have that:
\begin{eqnarray}\label{nierownosc2}
G(x,t+M)\leq C_0(1+|x+(t+M)|)^{-2}e^{-(\alpha-\epsilon)x-\theta_2^\epsilon(t+M)}.
\end{eqnarray}
Including $x=u-cu-f^{-1}(cu)$ and $t=f^{-1}(cu)$ into the equation \ref{lematZbyszkowy} and applying it to (\ref{numer2}) and (\ref{nierownosc2}) give:
\begin{eqnarray}\label{nierownosc3}
\lefteqn{\PP(\overline{\tau}(u-cu-f^{-1}(cu))\leq f^{-1}(cu))} \\ \nonumber&\leq&\frac{C_0(1+|(u-cu-f^{-1}(cu))+(f^{-1}(cu)+M)|)^{-2}}{h(0-,M)}\\
\nonumber&&\times e^{-(\alpha-\epsilon)(u-cu-f^{-1}(cu))-\theta_2^\epsilon(f^{-1}(cu)+M)}.
\end{eqnarray}
Sending $q\to\infty$ (note that then $h(0-,M)\to 1$) and then sending $M\to 0$
complete the proof.
\hfill $\square$
\\[1cm]

\begin{remark}\rm
For $\alpha\geq 0$ and $\beta>3$ one can take $\epsilon \downarrow 0$ in the assertion of Proposition \ref{upperbound}
to get upper bound independent of $\epsilon$ since then
$\theta_2^\epsilon\to\theta^0$ for some fixed $\theta^0$ and constant $C_0$ is bounded above.
\end{remark}

\texttt{Proof of Theorem \ref{LDP}.}
The assertion of the theorem follows from Propositions \ref{lowerbound} and \ref{upperbound} where we take $\epsilon \downarrow 0$ in latter one.
\hfill $\square$

\begin{remark}\rm
Let $f(t)=t^2$, $\beta >3$ and $\overline{F}(x)=e^{-\alpha x}$. In this case:
$$\lim_{u\to\infty}\frac{1}{u}\log \PP\left(\tau(u)\leq \sqrt{cu}\right)=-\alpha(1-c).$$
\end{remark}

\end{document}